\documentclass[12pt]{article}%
\usepackage{graphicx}
\usepackage{amsmath}
\usepackage{amsfonts}
\usepackage{amssymb}
\date{April, 2007}

\def\e{\epsilon}
\def\lf{\left}
\def\ri{\right}

\def\wt{\widetilde}

\def\p{\partial}
\newcommand\rr{{\mathbb R}}

\def\be{\begin{equation}}
\def\ee{\end{equation}}

\def\lf{\left}
\def\ri{\right}

\def\e{\epsilon}

\def\wt{\widetilde}

\def\p{\partial}

\def\wt{\widetilde}

\def\p{\partial}

\newcommand{\QED}{\hfill$\Box$\medskip}
\newtheorem{thm}{Theorem}[section]

\newtheorem{lem}{Lemma}[section]
\newtheorem{prop}{Proposition}[section]
\newtheorem{cor}{Corollary}[section]
\newtheorem{defn}{Definition}[section]
\newtheorem{rem}{Remark}
\numberwithin{equation}{section}
 \begin{document}
\title{Stability properties for the  higher dimensional catenoid in  $\rr^{n+1}$}
 \author{Luen-Fai Tam\thanks{Research partially supported by Earmarked Grant of Hong
Kong \#CUHK403005}\, and  Detang Zhou\thanks{Supported by CNPq of
Brazil.}}

\maketitle
\begin{abstract}  This paper concerns some stability properties
of higher dimensional catenoids  in $\rr^{n+1}$ with $n\ge 3$.  We
prove that higher dimensional catenoids have index one. We use
$\delta$-stablity for minimal hypersurfaces and show that the
catenoid is $\frac 2n$-stable and a  complete $\frac 2n$-stable
minimal hypersurface is a catenoid or a hyperplane provided the
second fundamental form satisfies some decay conditions.
\end{abstract}

\section{Introduction}

The catenoid in $\mathbb{R}^3$ is the only minimal surface of
revolution other than the plane. So it can be regarded as the
simplest minimal surface other than the plane. This motivates us to
study higher dimensional catenoids as complete minimal hypersurfaces
in higher dimensional Euclidean spaces $\rr^{n+1}$, $n\ge 3$. In
particular, we want to discuss  some stability properties of the
catenoids. Let us recall and introduce some notions of stability.

Let $M^n$ be a minimal hypersurface in $\rr^{n+1}$.  $M$ is said to
be {\it stable} if
\begin{equation}\label{stability1}
   \int_M\lf(|\nabla f|^2-|A|^2 f^2 \ri)\ge0
\end{equation}
for all $f\in C_0^\infty(M)$, where $|A|$ is the norm of the second
fundamental form of $M$. $M$ is said to be {\it weakly stable} if
(\ref{stability1}) is true for all $f\in C_0^\infty(M)$ with $\int_M
f=0$, see \cite{ccz}. Recall that in \cite{cm} it is defined that
$M$ is {\it $\delta$-stable} if
\begin{equation}\label{stability2}
   \int_M\lf(|\nabla f|^2-(1-\delta)|A|^2 f^2 \ri)\ge0
\end{equation}
for all $f\in C_0^\infty(M)$.

It is easy to see that $M$ is stable implies that $M$ is weakly
stable and $\delta$ stable when $\delta\ge 0$. If $M$ is a minimal
surface in $ \mathbb{R}^3$, then $M$ is always $\delta$-stable. By
\cite{fs}, $M$ is stable if and only if it there is a positive
solution $u$ of $(\Delta +|A|^2) u=0$. Hence $M$ is stable implies
that the universal cover of $M$ is stable. Similarly, $M$ is
$\delta$-stable implies that the universal cover of $M$ is also
$\delta$-stable. However, this is not true for weakly stable minimal
hypersurface.

In $\rr^3$, the catenoid is not stable by \cite{Ln}. In
$\rr^{n+1}(n\ge 3)$, it is proved in \cite{csz} that a complete
stable minimal hypersurface must have only one end. So  a catenoid
in  $\rr^{n+1}(n\ge 3)$ is not stable since it has two ends. In
fact, it is not even weakly stable (see \cite{ccz}). It is an
interesting question  to find the index of catenoids, which
measures the degree of instability. Using the Gauss map, it was
proved in \cite[p.131-132]{fc} that catenoids in $\mathbb{R}^3$
have index 1. It is known that a complete minimal surface in
$\mathbb{R}^3$ has finite index if and only if it has finite total
curvature, see \cite{fc}.   In \cite{Schoen1983}, Schoen proved
that the only complete nonflat embedded minimal surfaces in
$\mathbb{R}^3$ with finite total curvature    and with two ends
are the catenoids. It was also proved in \cite{LR} that only index
one complete minimal surfaces are the catenoid and Enneper surface
and catenoid is the only embedded minimal surface with index 1.
Although it has been believed   that a higher dimensional catenoid
also has index one, we have not found a reference for a proof. The
idea of using Gauss map in \cite{fc} does not work for higher
dimension. In this work, using a different method we prove that
the index of higher dimensional cateniod is indeed one, see
Theorem \ref{cateniodindex}. We would like to point out Choe
\cite{choe} has constructed higher dimensional Enneper's
hypersurfaces in $\rr^{n+1}$ when $n=3,4,5,6$. Different from the
catenoid, we don't even know whether it is of finite index.

It is well known that for $2\le n\le 6$, a complete area minimizing
hypersurface in $\rr^{n+1}$ must be a hyperplane. It is well known
by a result of do Carmo-Peng \cite{dcp79} and Fischer-Colbrie-Schoen
\cite{fs} independently that a complete stable minimal surface in
$\mathbb{R}^3$ is a plane.  On the other hand, for $3\le n\le 6$ it
is still an open question whether the condition of area minimizing
can be replaced by stability. In this direction, it was proved in
\cite{dcp80} (see also \cite{dcd}) that a complete stable minimal
hypersurface in $\rr^{n+1}$ is indeed a hyperplane under some
additional assumptions, for example: the norm of the second
fundamental form is square integrable. We will prove a similar
result  for catenoids which states that a complete $\frac 2n$-stable
minimal hypersurface in $\rr^{n+1}$ with $n\ge 3$ is a catenoid if
the norm of the second fundamental form satisfies certain decay
conditions. See Theorem \ref{quasistable}.
   As a corollary to this, we
show that if  a $\frac 2n$-stable complete proper immersed minimal
hypersurface $M^n$ in $\rr^{n+1}$ with $n\ge 3$  has least area
outside a compact set, and if the norm of the second fundamental
form is  square integrable  then $M$ is either a hyperplane or a
catenoid.

 The paper is  organized  as follows: in \S2, we introduce
 the definition and discuss some general properties of catenoids in higher
   dimensional Euclidean spaces. We will also prove that catenoids
   have index one. In \S3, we use the Simons' computation and the result in \cite{dcd}
    to give a characterization of catenoids. In \S4, we will
    discuss $\frac 2n$-stability and catenoids.

{\bf Acknowledgements.}  The second author would like to thank
 Bill Meeks, Harold Rosenberg and Rick Schoen for some
discussions during the 14th Brazilian geometry school in Salvador,
Brazil, 2006. We want to thank Rick Schoen for encouragements to
work for an index-one proof for the higher dimensional  catenoid.

\section{Catenoid and its index}\label{sectionindex}
In this section, we will  recall the definition of  catenoid, show
that it is $\frac 2n$-stable and compute its index.  Following do
Carmo and Dajczer \cite{dcd}, a catenoid is a complete rotation
minimal hypersurface in $\mathbb{R}^{n+1}$, $n\ge 2$ which is not a
hyperplane.  More precisely, let $\phi(s)$ be the solution of

\begin{equation}\label{catenoid1}
\left\{
  \begin{array}{ll}
    \frac{\phi''}{(1+{\phi'}^2)^\frac32}-
\frac{n-1}{\phi(1+{\phi'}^2)^\frac12}
  =0; \hbox{} \\
    \phi(0)  =\phi_0>0; \hbox{} \\
   \phi'(0)  =0. \hbox{}
  \end{array}
\right.
\end{equation}
 $\phi$ can be obtained as follows. Consider
\begin{equation}\label{catenoid2}
s=\int_{\phi_0}^\phi\frac{d\tau}{(a\tau^{2(n-1)}-1)^\frac12}
\end{equation}
where $a=\phi_0^{-2(n-1)}$. The integral in the right side of
(\ref{catenoid2}) is defined for all $\phi\ge \phi_0$. The function
$s(\phi)$ is increasing and if $n=2$,
 it maps $[\phi_0,\infty)$
onto $[0,\infty)$; if $n\ge3$, then it maps $[\phi_0, \infty)$ onto
$[0, S)$ where
$$
S=S(\phi_0)=\int_{\phi_0}^{+\infty}\frac{d\tau}{(a\tau^{2(n-1)}-1)^\frac12}<\infty.
$$
So $\phi(s)$ can be defined, and it is smooth up to $0$ such that
$\phi'(0)=0$. If we extend $\phi$ as an even function, then $\phi$
is smooth and satisfies (\ref{catenoid1}) on $\mathbb{R}$ in case
$n=2$ and on $(-S,S)$ in case $n\ge3$.

Let $\mathbb{S}^{n-1}$ be the standard unit sphere in
$\mathbb{R}^n$. A point $\omega\in \mathbb{S}^{n-1}$ can also be
considered as the unit vector $\omega$ in $\mathbb{R}^n$ which in
turn is identified as the hyperplane $x_{n+1}=0$ in
$\mathbb{R}^{n+1}$.

\begin{defn}\label{catenoiddef} A catenoid in
$\mathbb{R}^{n+1}$ is the hypersurface defined by the embedding:
$$
F:I\times\mathbb{S}^{n-1}\to \mathbb{R}^{n+1}
$$
with $F(s,\omega)=(\phi(s)\omega,s)$, where $I=\mathbb{R}$ if
$n=2$ and $I=(-S(\phi_0), S(\phi_0))$ if $n\ge3$, and $\phi$ is
the solution of (\ref{catenoid1}).
\end{defn}
 A hypersurface obtained by a rigid motion
of the hypersurface in the definition will also called a catenoid.
In case $n=2$, this is the standard catenoid in $\mathbb{R}^3$.
From now on, we are interested in the case that $n\ge3$.
\begin{prop}\label{basicproperties} Let $M$ be a catenoid
in $\mathbb{R}^{n+1}$ as in Definition \ref{catenoiddef}, $n\ge3$.
We have:
\begin{itemize}
  \item [(i)] $M$ is complete.
  \item [(ii)] The principal curvatures are
$\lambda_1=-\frac{\phi''}{(1+{\phi'}^2)^\frac32}$,
$\lambda_2=\cdots=\lambda_n=\frac{1}{\phi(1+{\phi'}^2)^\frac12}$.
  \item[(iii)] $M$ is minimal.
\item[(iv)] The norm $|A|$ of the second fundamental form $A$ of
$M$ is nowhere zero. Moreover, $|A|$ satisfies
\begin{equation}\label{catenoidsecond1}
  |A|\Delta |A|+|A|^4=\frac2n|\nabla|A||^2.
\end{equation}
\item[(v)] $M$ is symmetric with respect to the hyperplane
$x_{n+1}=0$ and is invariant under $O(n)$ which is the subgroup of
orthogonal transformations on $\mathbb{R}^{n+1}$ which fix the
$x_{n+1}$ axis. \item[(vi)] The part $\{x\in M|\ x_{n+1}\ge0\}$
and  the part $\{x\in M|\ x_{n+1}\le0\}$ are graphs over a subset
of $\{x_{n+1}=0\}$. \item[(vii)] Let $P$ be a hyperplane
containing the $x_{n+1}$ axis. Then $P$ divides $M$ into two
parts, each is a graph over $P$.
\end{itemize}
\end{prop}
\begin{proof}\  (i), (v), (vi) and (vii) are immediate consequences of
the definition. Let
$N=\frac{1}{(1+{\phi'}^2)^\frac12}(\omega,-\phi').$ Here and below
$'$ and $''$ are derivatives with respect to $s$. Then $N$ is the
unit normal of $M$. Let $D$ be the covariant derivative operator
in $\mathbb{R}^{n+1}$. Then
$$
D_{\frac{\p}{\p
s}}N=-\frac{\phi''}{(1+{\phi'}^2)^\frac32}(\phi'\omega,1)=-\frac{\phi''}{(1+{\phi'}^2)^\frac32}\frac{\p}{\p
s}.
$$
Suppose $(t_1,\dots,t_{n-1})$ are local coordinates of
$\mathbb{S}^{n-1}$, then
$$
D_{\frac{\p}{\p
t_i}}N=\frac{1}{(1+{\phi'}^2)^\frac12}(\frac{\p}{\p
t_i}\omega,0)=\frac{1}{\phi(1+{\phi'}^2)^\frac12}(\phi\frac{\p}{\p
t_i}\omega,0)=\frac{1}{\phi(1+{\phi'}^2)^\frac12}\frac{\p}{\p
t_i}.$$ From these (ii) follows.

(iii) follows from (ii) and (\ref{catenoid1}).

(iv) First note that (\ref{catenoid2}) implies
$\phi'=(a\phi^{2(n-1)}-1)^{\frac 12},$ then,
\begin{equation}\label{secondfundamentalform1}
  \begin{split}
       |A|^2&=\frac{n(n-1)}{\phi^2(1+{\phi'}^2)}\\
            & =n(n-1)\phi_0^{2(n-1)}\phi^{-2n}.
        \end{split}
\end{equation}
Hence $|A|>0$ everywhere because $\phi\ge\phi_0>0$. On the other
hand, the metric on $M$ in the coordinates $s,\omega$ is given by
\begin{equation}\label{metric}
g=(1+{\phi'}^2)ds^2+\phi^2 g_{\mathbb{S}^{n-1}}
\end{equation}
where $g_{\mathbb{S}^{n-1}}$ is the standard metric on
$\mathbb{S}^{n-1}$. Then
\begin{equation}\label{Laplacian}
\begin{split}
       \Delta \phi&=(1+{\phi'}^2)^{-1}\phi''+
\lf[(1+{\phi'}^2)^{-1}\ri]'\phi'+(1+{\phi'}^2)^{-1}\phi'
\lf[\log\lf((1+{\phi'}^2)^{\frac12}\phi^{n-1}\ri)\ri]'\\
            & =(1+{\phi'}^2)^{-1}\phi''-2(1+{\phi'}^2)^{-2}{\phi'}^2
\phi''+(1+{\phi'}^2)^{-1}\phi' \lf[\frac{\phi'\phi''}{
1+{\phi'}^2}+\frac{(n-1)\phi'}{\phi}\ri]\\
&=\frac{\phi''}{ (1+{\phi'}^2)^2}+(n-1)\frac{|\nabla
\phi|^2}\phi\\&=(n-1)\frac{1}{ \phi(1+{\phi'}^2)
}+(n-1)\frac{|\nabla \phi|^2}\phi
        \end{split}
\end{equation}
where $\nabla$ is the covariant derivative of $M$ and we have used
(\ref{catenoid1}). (\ref{catenoidsecond1}) follows from
(\ref{secondfundamentalform1}) and (\ref{Laplacian}) by a direct
computation.\QED
\end{proof}

By (iv) of the Proposition,  we see that
$$
\Delta |A|^{\frac n{n-2}}+\frac  {n-2}n|A|^2 |A|^{\frac n{n-2}}=0
$$
and $|A|^{\frac n{n-2}}>0$. Hence the catenoid is $\frac 2n$-stable
by \cite{fc}.

\begin{thm}\label{cateniodindex}Let $M$ be a catenoid in $\rr^{n+1}$. Then index of $M$ is 1.
\end{thm}

\begin{proof}\  It is well known that $M$ is not stable.  One can also  use the result of  Cao, Shen and
Zhu\cite{csz}. They proved that any complete stable minimal
hypersurfaces in $\rr^{n+1}$ has only one end, since the catenoid
has two ends, thus the index of $M$ is at least 1. We only need to
prove that its index is at most 1. Recall that the stability
operator is written as
\begin{equation*}
    L=\Delta+|A|^2.
\end{equation*}
For $M$ above $|A|^2(x)$ is  an even function depending only on
$r$. From the fact that $M$ is unstable it follows that
$\lambda_1(L)<0$. we now show that the second eigenvalue
$\lambda_2^D(L)\ge 0$ of $L$ on any bounded domain $D\subset M$.
Assume for the sake of contradiction that it were not true, we can
find a domain $D(R)=(-R,R)\times \mathbb{S}^{n-1}$ such that
$\lambda_2^{D(R)}(L)< 0$. Here $0<R<S$ and $S=S(\phi_0)$ is as in
Definition \ref{catenoiddef}.  That is to say that there is a
function $f$ satisfying
\begin{equation}\label{eq2}
    \left\{
  \begin{array}{ll}
    Lf=-\lambda_2f, & \hbox{ in $D(R)$;} \\
    f|_{\partial D(R)}=0. & \hbox{}
  \end{array}
\right.
\end{equation}
We claim that $f$ depends only on $r$. For any unit vector $v\in
\mathbb{S}^n$, and $v\perp (1,0,\cdots, 0)$, denote by $\pi_v$ the
hyperplane
\begin{equation*}
    \{p\in \rr^{n+1}, \langle p,v\rangle=0\}.
\end{equation*}
Let $\sigma_v$ be the reflection with respect to $\pi_v$. Define
function $\varphi_v(r,\theta)= f(r,\theta)-f_v(r,\theta)$ where
$f_v(p):=f(\sigma_v(p))$ for any $p\in D(R)$. Since
\begin{equation}\label{eq3}
    \Delta f(r,\theta)=\frac{\partial^2f}{\partial
r^2}+\frac{a'(r)}{a(r)}\frac{\partial f}{\partial
r}+\frac{1}{a^2(r)}\Delta_{S^n-1}f,
\end{equation}
 $f_v$ also satisfies (\ref{eq3}). Then
\begin{equation}\label{eq4}
    \left\{
  \begin{array}{ll}
    L\varphi_v=-\lambda_2\varphi_v, & \hbox{ in $D(R)$;} \\
    \varphi_v|_{D(R)\cap \pi_v}=0. & \hbox{}
  \end{array}
\right.
\end{equation}
Denote
\begin{equation*}
    D_v^+(R):=\{p\in \rr^{n+1}, \langle p,v\rangle>0\}.
\end{equation*}
Then $D_v^+(R)$ is a minimal graph over a domain in $\pi_v$ thus
is stable. From (\ref{eq4}) and $\lambda_2<0$, we conclude that
$\varphi_v\equiv 0$. It is a well-known fact that any element in
orthogonal group $O(n-1)$ can be expressed as a composition of
finite number of reflections, we know that $f$ is rotationally
symmetric.

Since $f$ is the second eigenfunction of $L$, it changes sign, so
there exists a number $r_0\in (-R,R)$ such that $f(r_0)=0.$ Assume
without loss of generality that $r_0\ge 0$. We take
$D(r_0,R):=\{p=(r,\theta)\in D(R), r\in (r_0, R)\}.$  Again $f$ is
an eigenfunction of $L$ on $D(r_0,R)$. Again we know that
$D(r_0,R)$ is a  minimal graph which contradicts the fact that
$\lambda_2<0$ because $f$ cannot be identically zero in
$D(r_0,R)$. The contradiction shows the index of $M$ is $1$.\QED
\end{proof}

\section{Simons' equation and catenoid}\label{sectionsimons}

By Proposition \ref{basicproperties}, the norm of the second
fundamental of a catenoid is nowhere zero and satisfies
(\ref{catenoidsecond1}).  In this section, we will prove that a
complete non flat minimal hypersurface in $\rr^{n+1}$ satisfying
(\ref{catenoidsecond1}) must be a catenoid.  Let us recall the
Simons' computation on the second fundamental form of a minimal
hypersurface in Euclidean space.

Let $M$ be an  $n$-dimensional manifold immersed in $\rr^{n+1}$. Let
$A$ be its second fundamental form and $\nabla A$ be its covariant
derivative. Let $h_{ij}$ and $h_{ijk}$ be the components of $A$ and
$\nabla A$ in an orthonormal frame.

By Proposition \ref{basicproperties}(iv), we see that Simons
inequality becomes equality for   catenoids. We will prove that the
converse is also true. We first prove a lemma.

\begin{lem}\label{basicequation}
Let $M$ be an immersed oriented minimal hypersurface in
$\rr^{n+1}$. At a point where the norm of the second fundamental
form $|A|>0$, we have
\begin{equation}\label{eqnsimons}
    |A|\Delta |A|+|A|^4=\frac 2n|\nabla|A||^2+E.
\end{equation}
with $E\ge0$. Moreover, in an orthonormal frame $e_i$ such that
$h_{ij}=\lambda_i\delta_{ij}$, then $E=E_1+E_2+E_3$, where
\begin{equation}\label{E's}
    \left\{
      \begin{array}{ll}
        E_1 &=\sum_{j\ne i,k\ne i, k\ne j}h_{ijk}^2 , \\
        E_2 &= \frac2n\sum_{j\ne i,k\ne i, k\ne j}(h_{kki}-h_{jji})^2, \\
        E_3 &=(1+\frac 2n)|A|^{-2} \sum_k\sum_{i\ne
    j}(h_{ii}h_{jjk}-h_{jj}h_{iik})^2.
      \end{array}
    \right.
\end{equation}
\end{lem}
\begin{proof}\ At a point $p$ where $|A|>0$, choose an
orthonormal frame such that $h_{ij}=\lambda_i\delta_{ij}$. Since
$M$ is minimal, then by \cite[(1.20),(1.27)]{ssy},     for $|A|>0$
we have:
\begin{equation}\label{Simons1}
   |A|\Delta |A|+|A|^4=\sum_{i,j,k=1}^n h_{ijk}^2-|\nabla|A||^2.
\end{equation}
Now,
\begin{equation}\label{Simons2}
\begin{split}
  |\nabla |A||^2 &=\left[ \sum_k(\sum_{i}h_{ii}h_{iik})^2\right] |A|^{-2} \\
   &=  \left[ \sum_k(\sum_{i}h_{ii}^2\sum_ih_{iik}^2)-\sum_k\sum_{i\ne j}(h_{ii}h_{jjk}-h_{jj}h_{iik})^2\right] |A|^{-2}\\
   &= \sum_{k,i}h_{iik}^2-\left[\sum_{i\ne j}(h_{ii}h_{jjk}-h_{jj}h_{iik})^2\right]
   |A|^{-2}
\end{split}
\end{equation}
where we have used the fact that $M$ is minimal. On the other
hand,
\begin{equation}\label{Simons3}
   \begin{split}
     \sum_{k,i}h_{iik}^2 & = \sum_{k\ne i}h_{iik}^2+\sum_ih_{iii}^2\\
       &= \sum_{k\ne i}h_{iik}^2+\sum_i(\sum_{j\ne
i}h_{jji})^2\\
&=\sum_{k\ne i}h_{iik}^2+\sum_i\left[(n-1)\sum_{j\ne i}h_{jji}^2-
\sum_{j\ne i,k\ne i, k\ne j}(h_{kki}-h_{jji})^2\right]\\
 &=n \sum_{k\ne i}h_{iik}^2-\sum_{j\ne i,k\ne i, k\ne j}(h_{kki}-h_{jji})^2
 \end{split}
\end{equation}
Combining this with (\ref{Simons2})
\begin{equation}\label{Simons4}
\sum_{k\ne i}h_{iik}^2=\frac1n\lf[|\nabla |A||^2 +\left(\sum_{i\ne
j}(h_{ii}h_{jjk}-h_{jj}h_{iik})^2\right)
   |A|^{-2}+\sum_{j\ne i,k\ne i, k\ne
j}(h_{kki}-h_{jji})^2\ri]
\end{equation}
 Note that since $\rr^{n+1}$ is flat, we have
$h_{ijk}=h_{ikj}$, (see \cite[(1.13)]{ssy} for example). By
(\ref{Simons2}),
\begin{equation}\label{Simons5}
    \begin{split}
\sum_{i,j,k}h_{ijk}^2& -|\nabla|A||^2\\
&=\sum_{j\ne i,k\ne i, k\ne j}h_{ijk}^2+\sum_{i\ne
k}h_{iik}^2+\sum_{i\ne k}h_{iki}^2+\sum_{i\ne
k}h_{ikk}^2+\sum_{i}h_{iii}^2-|\nabla|A||^2\\
&=\sum_{j\ne i,k\ne i, k\ne j}h_{ijk}^2+2\sum_{i\ne
k}h_{iik}^2+\sum_{i,k}h_{iik}^2-|\nabla|A||^2\\
&=\sum_{j\ne i,k\ne i, k\ne j}h_{ijk}^2+2\sum_{i\ne
k}h_{iik}^2+\left(\sum_{i\ne
j}(h_{ii}h_{jjk}-h_{jj}h_{iik})^2\right)
   |A|^{-2}
\end{split}
\end{equation}
(\ref{eqnsimons}) follows from (\ref{Simons1}), (\ref{Simons4})
and (\ref{Simons5}).
\end{proof}\QED

Since $E$ are nonnegative, we have the following Simons
inequality, see \cite{ssy}:
\begin{equation}\label{ineqsimons}
    |A|\Delta |A|+|A|^4\ge \frac 2n|\nabla|A||^2
\end{equation}
at the point $|A|>0$.

Now we are ready to prove the following:

\begin{thm}\label{thmcatenoid} Let $M^n(n\ge 3)$ be a non-flat  complete immersed minimal hypersurface in
$\rr^{n+1}$. If the   Simons inequality (\ref{ineqsimons}) holds
as an equation on all nonvanishing point of $|A|$ in $M$, then $M$
must be a catenoid.
\end{thm}
\begin{proof}\  Suppose $\Phi: M\to \rr^{n+1}$ is the minimal
immersion. Since $M$ is not a hyperplane, then $|A|$ is a
nonnegative continuous function which does not vanish identically.
Let $p$ be  a point such that $|A|(p)>0$. Then $|A|>0$ in a
connected open set $U$ containing $p$. Suppose that $|\nabla
|A||\equiv 0$ in $U$, then $|A|$ is a positive constant in $U$.
Since $|A|$ satisfies:
$$
|A|\Delta |A|+|A|^4=\frac 2n|\nabla|A||^2
$$
we have a contradiction. Hence there is a point in $U$ such that
$|\nabla |A||\ne0$. By shrinking $U$, we may assume that $|A|>0$
and $|\nabla|A||>0$ in $U$.
 By
(\ref{eqnsimons}) and the fact that (\ref{ineqsimons}) is an
equality in $U$, we conclude that $E=0$ in $U$.

Let $q\in U$. Choose an orthonormal frame at $q$ so that the
second fundamental form is diagonalized,
$h_{ij}=\lambda_i\delta_{ij}$. $E_2=0$ implies
\begin{equation*}
    h_{jji}=h_{kki}, \textrm{ for all }j\ne i,\quad k\ne i.
\end{equation*}
Combining with the minimal condition, we have
\begin{equation}\label{hiii}
    h_{iii}=-(n-1)h_{jji}, \textrm{ for all }j\ne i.
\end{equation}
Since $|\nabla |A||\ne 0$, then there exist $i_0$ and $j_0\ne i_0$
such that $h_{j_0j_0i_0}\ne 0$ hence $h_{i_0i_0i_0}\ne 0$. Suppose
for simplicity that $i_0=1$.

 $E_3=0$ implies
\begin{equation*}
    h_{ii}h_{jjk}=h_{jj}h_{iik},\textrm{ for all }i,j,k,
\end{equation*}
then
\begin{equation*}
    h_{11}h_{jj1}=h_{jj}h_{111}=-(n-1)h_{jj}h_{jj1},\textrm{ for all }j\ne
    1,
\end{equation*}
by (\ref{hiii}). So
\begin{equation}\label{eqnh11}
    h_{11}=-(n-1)h_{jj},\textrm{ for all }j\ne    1
\end{equation}
because $-(n-1)h_{jj1}=h_{111}\ne0$.  Hence the eigenvalues of
$h_{ij}$ are $\lambda$ with multiplicity $n-1$ and $-(n-1)\lambda$
with $\lambda\neq0$ because $|A|>0$. Hence in a neighborhood of of
$p$ the eigenvalues of $h_{ij}$ are of this form. By a result of
do Carmo and Dajczer \cite[Corollary 4.4]{dcd}, this neighborhood
is part of a catenoid. Hence $\Phi(M)$ is contained in a catenoid
$\mathcal{C}$ by minimality of the immersion. Since $M$ is
complete and $\Phi$ is a local isometry into the catenoid
$\mathcal{C}$ which is simply connected because $n\ge3$, $\Phi$
must be an embedding, see \cite[p.330]{Sp}. Hence $\Phi(M)$ is the
catenoid.\QED
\end{proof}

\section{$\frac 2n$-stability and catenoid}\label{sectionquasi}

In this section, we will prove that a complete immersed minimal
hypersurface in $\rr^{n+1}$, $n\ge3$ is a catenoid if it is $\frac
2n$-stable and if the second fundamental form satisfies some decay
conditions. We will also discuss the case when the minimal
hypersurface is area minimizing outside a compact set.

Following \cite{ssy}, let $M$ be a complete immersed minimal
 hypersurface in $\mathbb{R}^{n+1}$, $n\ge 3$. Assume there is a
 Lipschitz function $r(x)$ defined on $M$ such that
 $|\nabla r|\le 1$ a.e. Define $B(R)$ for $0<R<\infty$ by
$$
B(R)=\{x\in M|\ r(x)<R\}.
$$
Assume also that $B(R)$ is compact for all $R$ and
$M=\bigcup_{R>0}B(R)$. For example, $B(R)$ may be an intrinsic
geodesic ball or the intersection of an extrinsic ball with $M$.
In the later case, we assume that $M$ is proper.

\begin{thm}\label{quasistable} Let $M^n(n\ge 3)$ be a $\frac 2n$-stable complete immersed minimal hypersurface in
$\rr^{n+1}$. If
\begin{equation}\label{eqnlimit}
  \lim_{R\to +\infty} \frac
1{R^2}\int_{B(2R)\setminus B(R)}|A|^{\frac{2(n-2)}{n}}=0,
\end{equation}
then  $M$ is either a plane or a catenoid.
\end{thm}
\begin{proof}\  For any $\e>0$, let $u:=\lf(|A|^2+\e\ri)^{\frac\alpha2}$,
where $\alpha=\frac{n-2}{n}$. Then at the point $|A|>0$,
\begin{equation}\label{Bernstein1}
   \begin{split}
       \Delta u &=u\lf(\Delta \log u+|\nabla\log u|^2\ri)\\
&=\frac{\alpha
u}2\lf(\frac{\Delta|A|^2}{|A|^2+\e}-\frac{|\nabla|A|^2|^2}{(|A|^2+\e)^2}\ri)
+\frac{u\alpha^2}4\frac{|\nabla|A|^2|^2}{(|A|^2+\e)^2}\\
&=\alpha u\lf(\frac{\frac12\Delta
|A|^2}{|A|^2+\e}+(\alpha-2)\frac{|A|^2|\nabla|A||^2}{(|A|^2+\e)^2}\ri)\\
&=\alpha u\lf(\frac{(2-\alpha)|\nabla
|A||^2-|A|^4+E}{|A|^2+\e}+(\alpha-2)
\frac{|A|^2|\nabla|A||^2}{(|A|^2+\e)^2}\ri)\\
&\ge -\alpha u |A|^2 +\frac{\alpha u E}{ |A|^2+\e }
\end{split}
\end{equation}
where we have used (\ref{stability2})  and $E=E_1+E_2+E_3\ge0$. If
we extend $E$ to be zero for $|A|=0$, then it is easy to see that
the above inequality is still true.

 On the other hand, for any
function $\phi\in C_o^{\infty}(M)$,
\begin{eqnarray}\label{eqnineq1}
 \int_M \phi^2\frac{\alpha u E}{ |A|^2+\e }&\le& \int_M\phi^2u\left(\Delta u + \alpha|A|^2u \right) \nonumber \\
 &=&  -\int_M\phi^2|\nabla u|^2 -2\int_M\phi u\langle    \nabla u,\nabla \phi\rangle +\int_M\alpha|A|^2\phi^2 u^2 \nonumber \\
   &\le&  -2\int_M\phi u\langle    \nabla u,\nabla \phi\rangle  -\int_M\phi^2|\nabla u|^2 +\int_M|\nabla(\phi u)|^2 \nonumber \\
   &=& \int_M|\nabla\phi|^2u^2.
\end{eqnarray}
Here we have used (\ref{stability2}). Let $\phi$ be a smooth
function on $[0,\infty)$ such that $\phi\ge0$, $\phi=1$ on $[0,R]$
and $\phi=0$ in $[2R,\infty)$ with $|\phi'|\le \frac2R$. Then
consider $\phi\circ r$, where $r$ is the function in the definition
of $B(R)$.
\begin{equation}\label{ineqr}
 \int_{B(R)} \phi^2\frac{\alpha u E}{ |A|^2+\e } \le
 \int_{B(R)}\phi^2u\left(\Delta u + \alpha|A|^2u \right)\le
\frac 4{R^2}\int_{B(2R)\setminus B(R)}||A|^{\frac{2(n-2)}{n}}.
\end{equation}
Let $\e\to 0$ and then let  $R\to +\infty$, we conclude that $E=0$
whenever $|A|>0$. Thus the Simons' inequality becomes equality on
$|A|>0$.  By Theorem \ref{thmcatenoid}, we know that it is a
catenoid.
\end{proof}
\begin{rem} It should be remarked that  (\ref{eqnlimit})  is
satisfied when $M$ is a catenoid. In fact, using notation in the
Definition \ref{catenoiddef}, the metric is of the form
$$
g=(1+\phi'^2)ds^2+\phi^2 g_{\mathbb{S}^{n-1}}.
$$
Hence the distance function is of order $\phi$.  By
(\ref{secondfundamentalform1}), $|A|$ is  of order $\phi^{-n}$.
The volume of geodesic ball of radius $r\sim \phi$ is of order
$\phi^n$. From this it is easy to see that (\ref{eqnlimit}) is
true for $n\ge3$.
\end{rem}

We say that $M$ has least area outside a compact set (see
\cite{ssy}, p. 283]), if (i) $M$ is proper; and (ii)
 $M$ is the boundary of some open set $U$ in
$\mathbb{R}^{n+1}$ and there is $R_0>0$ such that for any open set
$\mathcal{O}$ in $\mathbb{R}^{n+1}$ with $\mathcal{O}\cap \wt
B(R_0)=\emptyset$ we have $|\p U\cap \mathcal{O}|\le |\p
\mathcal{O}\cap U|$. Here $\wt B(R_0)$ is the extrinsic ball in
$\mathbb{R}^{n+1}$ with center at the origin. If this is true,
then $M$ is stable outside a compact set and if $r$ is the
extrinsic distance, then
$$|B(4R)\setminus B(\frac12 R)|\le |\p\wt B(4R)|+|\p\wt
B(\frac12 R)|\le CR^n$$ if $R$ is large.
\begin{cor}\label{leastarea} Let $M^n$, $n\ge 3$ be a $\frac 2n$-stable
complete proper immersed minimal hypersurface in $\rr^{n+1}$. If
$M$ has least area outside a compact set and
\begin{equation}\label{integralcondition}
   \int_M|A|^p<\infty,
\end{equation}
for some $\frac{2(n-2)}n\le p\le 2$ then $M$ is either a plane or
a catenoid.
\end{cor}
\begin{proof}\ Suppose $|A|$ satisfies (\ref{integralcondition}).
Since $\frac{2(n-2)}n\le p\le 2$, we have:
\begin{equation}
   \begin{split}
 R^{-2}\int_{B(2R)\setminus B(R)}|A|^{\frac{2(n-2)}n}&\le
 CR^{-2}\lf(\int_M|A|^p\ri)^{\frac{2(n-2)}{pn}}R^{n-\frac{2(n-2)}p}\\
 &=C\lf(\int_M|A|^p\ri)^{\frac{2(n-2)}{pn}}R^{(n-2)(1-\frac
 2p)}\to 0
\end{split}
\end{equation}
as $R\to \infty$. The result follows from Theorem
\ref{quasistable}.
 \QED
\end{proof}

By \cite{Schoen1983}, the only nonflat complete minimal immersions
of $M^n\subset \rr^{n+1}$, which are regular at infinity and have
two ends, are the catenoids. By the corollary, we have the
following:

\begin{cor}\label{regular} A nonflat complete minimal immersion  of $M^n\subset
\rr^{n+1}$ with $n\ge 3$, which are regular at infinity and has more
than two ends, is not $\frac 2n$-stable.
\end{cor}

\bigskip\bigskip\noindent Luen-fai Tam\newline Department of
Mathematics\newline The Chinese University of Hong Kong\newline Shatin,
Hong Kong, China\newline e-mail: lftam@math.cuhk.edu.hk\\

\noindent Detang Zhou\newline Insitituto de Matematica\newline
Universidade Federal Fluminense- UFF\newline Centro, Niter\'{o}i,
RJ 24020-140, Brazil
\newline email: zhou@impa.br

 \end{document}